\documentclass [12pt,a4paper,leqno]{amsart}
\input amssymb.sty
\usepackage{exscale}
\usepackage{fancyhdr}
\usepackage{amsfonts}
\usepackage{amssymb}
\usepackage{ulem}

\newtheorem{thm}{Theorem}
\newtheorem{definition}{Definition}
\newtheorem{lem}[thm]{Lemma}
\newtheorem{corollary}[thm]{Corollary}
\newtheorem{proposition}[thm]{Proposition}

\theoremstyle{remark}

\newtheorem*{acknowledgment*} {Acknowledgment}

\def\Ln{\operatorname{Ln}}

\newenvironment{proofof}[1]
        {\pagebreak[2] \vspace{-1pt}{\bf Proof of #1.}  }
        {\hfill $\blacksquare$ \vspace{2pt}}

\renewcommand{\rho}{\varrho}
\renewcommand{\phi}{\varphi}
\renewcommand{\epsilon}{\varepsilon}

\def\beq{\begin{equation}}
\def\eeq{\end{equation}}
\def\beqar{\begin{eqnarray}}
\def\eeqar{\end{eqnarray}}
\def\beqaro{\begin{eqnarray*}}
\def\eeqaro{\end{eqnarray*}}
\def\bsat{\begin{thm}}
\def\esat{\end{thm}}
\def\blem{\begin{lem}}
\def\elem{\end{lem}}
\def\bkor{\begin{corollary}}
\def\ekor{\end{corollary}}
\def\bprop{\begin{proposition}}
\def\eprop{\end{proposition}}
\def\bdefin{\begin{definition}}
\def\edefin{\end{definition}}
\def\bbew{\begin{proof}}
\def\ebew{\end{proof}}
\def\bbewo{\begin{proofof}}
\def\ebewo{\end{proofof}}

\renewcommand{\rho}{\varrho}
\renewcommand{\phi}{\varphi}

\begin{document}

\thispagestyle{plain}

\title[The Natural Logarithm of a Natural Number]
{A New Formula for the Natural Logarithm of a Natural Number}

 \author{Shahar Nevo}
 \address{Bar-Ilan University,
Department of Mathematics, Ramat-Gan 52900, Israel}
 \email{nevosh@macs.biu.ac.il}

\begin{abstract}
For every natural number $T,$ we write $\Ln T$ as a series,
generalizing the known series for $\Ln 2.$
\end{abstract}



\thanks{This research is part of the European Science Foundation Networking
Programme HCAA and was supported by  Israel Science Foundation
Grant 395/07.}

\maketitle

\section{Introduction}

 The Euler-Mascheroni constant $\gamma$, \cite{1}, is given by
 the limit
 \begin{equation}\label{1}
 \gamma=\lim_{n\to\infty} A_n,
 \end{equation}
 where for every $n\ge 1,$ $A_n:= 1+\frac{1}{2}+\dots
 +\frac{1}{n}-\Ln n.$ An elementary way to show the convergence of
 $\{A_n\}_{n=1}^\infty$ is to consider the series
 $\sum_{n=0}^\infty(A_{n+1}-A_n)$. (Here $A_0:=0.)$ Indeed, by
 Lagrange's Mean Value Theorem, there exists for every $n\ge1$ a
 number $\theta_n,$ $0< \theta_n<1$ such that
 $$A_{n+1}-A_n=\frac{1}{n+1}-\Ln(n+1)+\Ln n=\frac{1}{n+1}-\frac{1}{n+\theta_n}=\frac{\theta_n-1}
 {(n+1)(n+\theta_n)},$$
and thus $0>A_{n+1}-A_n> \frac{-1}{n(n+1)}$ and the series
converges to some limit $\gamma.$

\section{The new formula}

Let $T\ge 2$ be an integer. We have
\begin{equation}\label{2}
A_{nT}=\sum_{k=0}^{n-1}\sum_{j=1}^T
\frac{1}{kT+j}-\Ln(nT)\underset{ n\to\infty}\rightarrow
\gamma.\end{equation} By subtracting \eqref{1} from \eqref{2} and
using $\Ln(nT)=\Ln n+\Ln T,$ we get
$$\sum_{k=0}^{n-1}\left(\sum_{j=1}^T\frac{1}{kT+j}-\frac{1}{k+1}\right)\underset
{n\to\infty}\rightarrow \Ln T,$$ that is,
\begin{equation}\label{3}
\Ln
T=\sum_{k=0}^\infty\left(\frac{1}{kT+1}+\frac{1}{kT+2}+\dots+\frac{1}{kT+(T-1)}-\frac{(T-1)}{kT+T}
\right).\end{equation} We observe that \eqref{3} generalizes the
formula $\Ln 2=1-\frac{1}{2}+\frac{1}{3}-\frac{1}{4}+\dots\, .$

We can write \eqref{3} also as
\begin{equation}\label{4}
\Ln
T=\left(1+\frac{1}{2}+\frac{1}{3}+\dots+\frac{1}{T}-1\right)+\left(\frac{1}{T+1}+\frac{1}{T+2}
+\dots+\frac{1}{2T}-\frac{1}{2}\right)+\dots\end{equation} and
this gives $\Ln T$ as a rearrangement of the conditionally
convergent series
$1-1+\frac{1}{2}-\frac{1}{2}+\frac{1}{3}-\frac{1}{3}+\dots.$ The
formula \eqref{4} holds also for $T=1.$ Formulas \eqref{3} and
\eqref{4} can be applied also to introduce $\Ln Q$ as a series for
any positive rational   $Q=\frac{M}{L}$ since $\Ln \frac{M}{L}=\Ln
M-\Ln L.$

Now, for any $k\ge0,$ the nominators of the $k$-th element in
\eqref{3} are the same and their sum is 0. This fact is not
random. For every constant $a_1,a_2,\dots a_T$, the sum
\begin{equation}\label{5}
S_T(a_1,\dots,a_T):=\sum_{k=0}^\infty\left(\frac{a_1}{kT+1}+\frac{a_2}{kT+2}+\dots+\frac{a_T}{(k+1)T}
\right)\end{equation} converges if and only $a_1+a_2+\dots
+a_T=0.$ This follows by comparison to the series
$\sum_{k=1}^\infty \frac{1}{k^2}<\infty.$  By \eqref{3} and the
notation \eqref{5}, $\Ln T=S_T(1,1,\dots,1,T-1).$

For $T\ge 2,$ let us denote by $\Sigma(T)$ the collection of all
sums of rational series of type \eqref{5}, i.e.,
$$\Sigma(T)=\big\{S_T(a_1,\dots,a_T): a_i\in  Q,\, 1\le i\le
T,a_1+\dots +a_T=0\big\}.$$ The collection $\Sigma(T)$ is a linear
space of real numbers over $\mathbb Q$ (or over the field of
algebraic numbers if we would define $\Sigma(T)$ to be with
algebraic coefficients instead of rational coefficients), and
$\dim\Sigma(T)\le T-1.$ A spanning set of $T-1$ elements of
$\Sigma(T)$ is
$$\big\{S_T(1,-1,0,0,\dots,0),S_T(0,1,-1,0,0,\dots,0),\dots,
S_T(0,\dots,0,1,-1)\big\}.$$ Also, if $T$ is not a prime number,
then $\dim\Sigma(T)<T-1.$ If $  Q=\frac{M}{L}$ is a positive
rational number and $P_1,P_2,\dots,P_k$ are all the prime factors
of $M$ and $L$ together, then $\Ln   Q \in \Sigma(P_1P_2\dots
P_k).$

We can get a non-trivial series for $x=0$: $\Ln 4=2\Ln
2=S_2(2,-2)=S_4(2,-2,2,-2)$, and also $\Ln (4)=S_4(1,1,1,-3).$
Hence
\begin{align*}0&=S_4(2,-2,2,-2)-S_4(1,1,1-3)=S_4(1,-3,1,1)\\
&=\left(\frac{1}{1}-\frac{3}{2}+\frac{1}{3}+\frac{1}{4}
\right)+\left(\frac{1}{5}-\frac{3}{6}+\frac{1}{7}+\frac{1}{8}\right)+\dots\,
.\end{align*}

\section{The integral approach}

The formula \eqref{3} can as well be deduced in the following way.
\begin{align}
\Ln T&=\lim_{x\to 1^-} \Ln(1+x+\dots +x^{T-1})=\lim_{x\to
1^-}\Ln\left(\frac{1-x^T}{1-x}\right)\nonumber\\
&=\lim_{x\to 1^-}(\Ln(1-x^T)-\Ln(1-x))=\lim_{x\to
1^-}\int_0^x\left[\frac{Tu^{T-1}}{u^T-1}+\frac{1}{1-u}\right]du\nonumber\\
&=\lim_{x\to 1^-}\int_0^x\frac{Tu^{T-1}-(1+u+\dots
+u^{T-1})}{u^T-1}du\nonumber\\
&=\lim_{x\to1^-}\int_0^x\frac{-1-u-u^2-\dots-u^{T-2}+(T-1)u^{T-1}}{u^T-1}du\nonumber\\
&=\lim_{x\to
1^-}\bigg[1\cdot\int_0^x\frac{u-1}{u^T-1}du+2\cdot\int_0^x
\frac{u^2-u}{u^T-1}du+3\int_0^x\frac{u^3-u^2}{u^T-1}du+\dots\nonumber\\
&\quad
+(T-2)\int_0^x\frac{u^{T-2}-u^{T-3}}{u^T-1}du+(T-1)\int_0^x\frac{u^{T-1}-u^{T-2}}{u^T-1}du\bigg]
\label{6}.
\end{align}

For every $1\le j\le T-1,$
\begin{align*}
\lim_{x\to 1^-}&\int_0^x\frac{u^j-u^{j-1}}{u^T-1}du =\lim_{x\to
1^-}\int_0^x\bigg[u^{j-1}\sum_{k=0}^\infty
u^{kT}-u^j\sum_{k=0}^\infty u^{kT}\bigg]du\\
&=\lim_{x\to 1^-}\int_0^x\bigg(\sum_{k=0}^\infty
u^{kT+j-1}-\sum_{k=0}^\infty u^{kT+j}\bigg)du =\lim_{x\to
1^-}\sum_{k=0}^\infty
\bigg(\frac{x^{kT+j}}{kT+j}-\frac{x^{kT+j+1}}{kT+j+1}\bigg).
\end{align*}
The series in the last expression converges at $x=1,$ and thus it
defines a continuous function in $[0,1]$ and so the limit is
\begin{equation}\label{7}
\int_0^1\frac{u^j-u^{j-1}}{u^T-1}du=\sum_{k=0}^\infty\left(\frac{1}{kT+j}-\frac{1}{kT+j+1}\right).
\end{equation}
By \eqref{6}, we now get that
\begin{align*}
\Ln
T&=\sum_{k=0}^\infty\left(\frac{1}{kT+1}-\frac{1}{kT+2}\right)+2\sum_{k=0}^\infty
\left(\frac{1}{kT+2}-\frac{1}{kT+3}\right)+\dots \\
&\quad+
(T-2)\sum_{k=0}^\infty\left(\frac{1}{kT+T-1}-\frac{1}{kT+T-1}\right)+(T-1)\sum_{k=0}^\infty\left(
\frac{1}{kT+T-1}-\frac{1}{(k+1)T}\right)\\
&=\sum_{k=0}^\infty\left(\frac{1}{kT+1}+\frac{1}{kT+2}+\dots+\frac{1}{kT+T-1}-\frac{(T-1)}{(k+1)T}
\right),
\end{align*}
and this is formula \eqref{3}.

If we put $T=3,$ $j=1$ into \eqref{7}, we get that
\begin{equation}\label{8}
\int_0^1\frac{u-1}{u^3-1}du=\sum_{k=0}^\infty\left(\frac{1}{3k+1}-\frac{1}{3k+2}\right)=S_3(1,-1,0).
\end{equation}
On the other hand,
$$\int\frac{u-1}{u^3-1}du=\frac{2}{\sqrt
3}\arctan\left(\frac{2u+1}{\sqrt 3}\right),$$ and together with
\eqref{8}, this gives
$$\frac{2}{\sqrt{3}}\left(\arctan\sqrt{3}-\arctan
\frac{1}{\sqrt{3}}\right)=S_3(1,-1,0)$$ or
$$\pi=3\sqrt {3}\cdot
S_3(1,-1,0)=3\sqrt{3}\left[\left(\frac{1}{1}-\frac{1}{2}\right)+\left(\frac{1}{4}-\frac{1}{5}\right)
+\left(\frac{1}{7}-\frac{1}{8}\right)+\dots\right].$$

\bibliographystyle{amsplain}

\end{document}